\documentclass[11pt,bezier]{article}
\usepackage{amsmath,mathrsfs,graphicx,amssymb,amsfonts,color}

\font\bg=cmbx10 scaled\magstep1

\font\small=cmr8

\newtheorem{newlemma}{{\bf Lemma}}

\newenvironment{lema}{\begin{newlemma}{\hspace{-0.5
em}{\bf.}}}{\end{newlemma}}

\newtheorem{newteorem}{{\bf Theorem}}

\newenvironment{teorem}{\begin{newteorem}{\hspace{-0.5
em}{\bf.}}}{\end{newteorem}}

\newtheorem{newkorolari}{{\bf Corollary}}

\newenvironment{korolari}{\begin{newkorolari}{\hspace{-0.5
em}{\bf.}}}{\end{newkorolari}}

\newtheorem{newdefine}{{\bf Definition}}

\newtheorem{newquestion}{{\bf Question}}

\newtheorem{newkonjek}{{\bf Conjecture}}

\newtheorem{newexample}{{\bf Example}}

\begin{document}
\tolerance=10000
\baselineskip18truept
\newbox\thebox
\global\setbox\thebox=\vbox to 0.2truecm{\hsize
0.15truecm\noindent\hfill}
\def\boxit#1{\vbox{\hrule\hbox{\vrule\kern0pt
\vbox{\kern0pt#1\kern0pt}\kern0pt\vrule}\hrule}}
\def\qed{\lower0.1cm\hbox{\noindent \boxit{\copy\thebox}}\bigskip}
\def\ss{\smallskip}
\def\ms{\medskip}
\def\bs{\bigskip}
\def\c{\centerline}
\def\nt{\noindent}
\def\ul{\underline}
\def\ol{\overline}
\def\lc{\lceil}
\def\rc{\rceil}
\def\lf{\lfloor}
\def\rf{\rfloor}
\def\ov{\over}
\def\t{\tau}
\def\th{\theta}
\def\k{\kappa}
\def\l{\lambda}
\def\L{\Lambda}
\def\g{\gamma}
\def\d{\delta}
\def\D{\Delta}
\def\e{\epsilon}
\def\lg{\langle}
\def\rg{\rangle}
\def\p{\prime}
\def\sg{\sigma}
\def\ch{\choose}

\newcommand{\ben}{\begin{enumerate}}
\newcommand{\een}{\end{enumerate}}
\newcommand{\bit}{\begin{itemize}}
\newcommand{\eit}{\end{itemize}}
\newcommand{\bea}{\begin{eqnarray*}}
\newcommand{\eea}{\end{eqnarray*}}
\newcommand{\bear}{\begin{eqnarray}}
\newcommand{\eear}{\end{eqnarray}}

\centerline{\Large \bf  On the edge cover polynomial of certain graphs}
\vspace{.3cm}

\bigskip

\baselineskip12truept
\centerline{ S. Alikhani$^{}${}\footnote{\baselineskip12truept\it\small
Corresponding author. E-mail: alikhani@yazd.ac.ir} and S. Jahari}
\baselineskip20truept
\centerline{\it Department of Mathematics, Yazd University}
\vskip-8truept
\centerline{\it  89195-741, Yazd, Iran}
\baselineskip20truept

\vskip-0.2truecm
\nt\rule{16cm}{0.1mm}

\nt{\bg ABSTRACT}
\medskip

\baselineskip14truept

\nt{ Let $G$ be a simple graph of order $n$ and size $m$.
The edge covering of  $G$ is a set of edges such that every vertex of $G$ is incident to at least one edge of the set.
 The edge cover polynomial of $G$ is the polynomial
$E(G,x)=\sum_{i=\rho(G)}^{m} e(G,i) x^{i}$,
where $e(G,i)$ is the number of edge coverings of $G$ of size $i$, and
$\rho(G)$ is the edge covering number of $G$. In this paper we study the
edge cover polynomials of cubic graphs of order $10$.
We show that all cubic graphs of order $10$ (especially the Petersen graph) are
determined uniquely by their edge cover polynomials. Also we construct   infinite families of graphs whose
edge cover polynomials have only roots $-1$ and $0$. }

\ms

\nt{\bf Mathematics Subject Classification:} {\small 05C60.}
\\
{\bf Keywords:} {\small Edge cover polynomial; edge covering; equivalence class; cubic graph; corona.}

\nt\rule{16cm}{0.1mm}

\baselineskip20truept

\section{Introduction}

\nt Let $G=(V,E)$ be a simple graph. The {\it order} and the size of $G$ is the number of vertices
and the number of edges of $G$, respectively.
For every graph $G$ with no isolated vertex, an edge covering of $G$ is a set of edges of $G$ such that every
vertex is incident with at least one edge of the set. In other words, an edge covering of a graph is a set
of edges which together meet all vertices of the graph. A minimum edge covering is an edge covering
of the smallest possible size. The edge covering number of $G$ is the size of a minimum edge covering of
$G$ and denoted by $\rho(G)$. In this paper we let $\rho(G) = 0$, if $G$ has some isolated vertices.
For a detailed treatment of these parameters, the reader is referred to~\cite{saeid1,Bond,GRo}.
Let $\mathscr{E}(G,i)$ be the family of all edge coverings of a graph $G$ with cardinality $i$ and
let $e(G,i)=|{\mathscr{E}}(G,i)|$.
The {\it edge cover polynomial} $E(G,x)$ of $G$ is defined as
\[
E(G, x)=\sum_{ i=\rho(G)}^{m} e(G, i) x^{i},
\]
where $\rho(G)$ is the edge covering number of $G$. Also, for a graph $G$ with some isolated vertices
we define $E(G, x) = 0$. Let $E(G, x) = 1$, when both order and size of $G$ are zero (see \cite{saeid1}).

\nt For two graphs $G$ and $H$, the corona $G\circ H$ is the graph arising from the
disjoint union of $G$ with $| V(G) |$ copies of $H$, by adding edges between
the $i$th vertex of $G$ and all vertices of $i$th copy of $H$ \cite{Fruc}. The corona $G\circ K_1$, in particular, is the graph constructed from a copy of $G$, where for each vertex $v\in V(G)$, a new vertex $u$ and a pendant edge $\{v, u\}$ are added.
It is easy to see that the corona operation of two graphs does not have the commutative property. 

\nt Two graphs $G$ and $H$ are said to be {\it edge covering equivalent},
or simply ${\mathscr{E}}$-equivalent, written $G\sim H$, if
$E(G, x)=E(H, x)$. It is evident that the relation $\sim$ of being
${\mathscr{E}}$-equivalence
is an equivalence relation on the family ${\cal G}$ of graphs, and thus ${\cal G}$ is partitioned into equivalence classes,
called the {\it ${\mathscr{E}}$-equivalence classes}. Given $G\in {\cal G}$, let
\[
[G]=\{H\in {\cal G}:H\sim G\}.
\]
We call $[G]$ the equivalence class determined by $G$.
A graph $G$ is said to be {\it edge covering unique}, or simply {\it ${\mathscr{E}}$-unique}, if $[G]=\{G\}$.

\nt For every vertex $v \in V(G)$, the degree of $v$ is the number of edges incident
with $v$ and is denoted by $d_G(v)$. For simplicity we write $d(v)$ instead of $d_G(v)$. Let $a_i(G)$ be the
number of vertices of $G$ with degree $i$.
The minimum degree of $G$ is denoted by
$\delta(G)$. A graph $G$ is called {\it $k$-regular} if all
vertices have the same degree $k$.
One of the famous graphs is
the Petersen graph which is a symmetric non-planar cubic graph. In the study of edge cover
polynomials, it is interesting to investigate the edge coverings and edge cover polynomial of
this graph. We denote the Petersen graph by $P$.

\nt In \cite{saeid1} authors has characterized all graphs whose edge cover polynomials have exactly
one or two distinct roots and moreover they proved that these roots
are contained in the set $\{-3,-2,-1, 0\}$. In this paper, we construct some  infinite families of graphs whose
 edge cover polynomials have only roots $-1$ and $0$. Also  similar to \cite{turk}, we study the edge coverings and edge cover
polynomials of cubic graphs of order $10$. As a consequence, we show that the all cubic graphs of order $10$ (especially the
 Petersen graph) are
determined uniquely by their edge cover polynomials.

\nt  In the next section, we construct some  infinite families of graphs whose edge cover polynomials have only roots $-1$ and $0$.
In Section 3, we obtain the edge cover polynomial of the Petersen graph. In section $4$, we list all $\rho$-sets of connected cubic graphs of order $10$. This list will be used to study the ${\mathscr{E}}$-equivalence of these graphs
in the last section. In the last section, using a procedure we study the ${\mathscr{E}}$-equivalence classes
of all cubic graphs of order $10$ and we show
that all cubic graphs of order $10$ are determined uniquely by their  edge cover polynomial.

\section{Families of graphs with edge cover root $-1$ and $0$}

\nt In this section we construct  infinite families of graphs whose  edge cover polynomials have only roots $-1$ and $0$.

\nt Let $G$ be any graph with vertex set $\{v_1, \cdots , v_n\}$. Add $n$ new vertices $\{u_1, \cdots , u_n\}$ and join $u_i$ to $v_i$ for every $i, 1 \leq i \leq n$. By definition we denote this new graph by $G\circ K_1$.  We start this section with the following Lemma.

\begin{lema}\label{lemma4}
For any graph $G$ of order $n$, $\rho(G\circ K_1)=n$.
\end{lema}
\nt{\bf Proof.} If $E$ is an edge covering of $G\circ K_1$, then for every $i$, $\{v_i,u_i\} \in E$. This implies that
$| E |\geq n$. Since $\{\{v_1, u_1\}, \cdots , \{v_n, u_n\}\}$ is an edge covering of $G\circ K_1$, we have the result.\quad\qed

\nt By Lemma \ref{lemma4}, $e(G\circ K_1, i)= 0$, for every $i, i < n$, so we shall compute $e(G\circ K_1, i)$ for each $i,
n \leq i \leq n+| E(G) |$.

\begin{teorem}
For any graph $G$ of order $n$ and size $m$. For every $i, n \leq i \leq n+m$, we have $e(G\circ K_1, i)={m \choose i-n}$.
Hence $E(G\circ K_1, x) = x^n(1+x)^m$.
\end{teorem}
\nt{\bf Proof.}
Suppose that $E$ is an edge covering of $G\circ K_1$ of size i, and $| E \cap E(G) |= j$, for $0\leq j\leq m$. Also $E$ contains these edges $\{v_1, u_1\}, \cdots , \{v_n, u_n\}.$ For addition $j$ edges of $E(G)$ to an edge covering $E$ of size $i$, we have ${m \choose j}$ possibilities.\quad\qed

\nt We have the following corollary:

\begin{korolari} \label{twor}
Let $i\in \mathbb{N}$.  For any graph $G$ of order $n$ and size $m$, $E(G\circ \overline{{K_i}}, x) = x^{in}(1+x)^m$.

\end{korolari}

 \nt In \cite{saeid1} authors has characterized all graphs whose edge cover polynomials have exactly
one or two distinct roots and moreover they proved that these roots
are contained in the set $\{-3,-2,-1, 0\}$. The following corollary is an immediate consequence of Corollary \ref{twor}.

\begin{korolari}
The edge cover polynomial of every graph $H$ in
the family $\{G\circ \overline{K_i}, (G\circ \overline{K_i})\circ \overline{K_i}, ((G\circ \overline{K_i})\circ \overline{K_i})\circ \overline{K_i},\cdots \}$ have only two roots $-1$ and $0$.
\end{korolari}

\nt As we know the edge cover  polynomials of binary graph operations, such as  join, corona, Cartesian and lexicographic product
 have not been studied. Let us to state the following problem for future studies.

\nt {\bf Problem.} Compute and present  recurrence formulae
and properties of the edge cover  polynomials of families of graphs obtained by
various products.

\section{Edge cover polynomial of the Petersen graph}

\nt In this section we shall investigate the edge cover polynomial of the Petersen graph.
First, we state some properties of the edge cover polynomial of a graph.

\begin{lema}\label{lemma1} {\rm \cite{saeid1}}
Let $G$ be a graph of order $n$ and size $m$ with no isolated vertex. If
$E(G, x)=\sum_{ i=\rho(G)}^{m} e(G,i) x^{i}$, then the following hold:
\begin{enumerate}
\item[(i)] $E(G, x)$ is a monic polynomial of degree $m$.
\item[(ii)] $n\leq 2\rho(G)$
\item[(iii)] For $i = m - \delta + 1, . . . ,m$, $e(G, i) ={m \choose i}$. Moreover, if $ i_0=\min\{i | e(G,i)={m\choose i}\},$ then $\delta= m-i_0+1.$
\item[(iv)] If $G$ has no connected component isomorphic to $K_2$, then $a_{\delta}(G) ={m \choose {m-\delta}} - e(G,m - \delta)$.
\end{enumerate}
\end{lema}

\begin{korolari}\label{korolari1}
Let $G$ and $H$ be two graphs (with no isolated vertex) of size $m_1$ and $m_2$, respectively. If
$E(G, x) = E(H, x)$, then $\rho(G) =\rho(H), m_1 = m_2$ and $\delta(G) = \delta(H).$
\end{korolari}

\nt The next theorem shows that the edge cover polynomial of a graph determines the regularity of graph.
\begin{teorem}\label{theorem1}{\rm \cite{saeid1}}
Let $G$ be a graph. If $\delta\geq 2$, then $G$ is regular if and only if
\[
e(G,m-\delta)={m\choose \delta}-\frac{2m}{\delta},
\]
where $m$ is the size of $G$.
\end{teorem}

\nt The following theorem states that if $G$ is a regular graph, then every element of $[G]$ is a regular graph too.

\begin{teorem}\label{theorem2}{\rm \cite{saeid1}}
Let $G$ be a $k$-regular graph of order $n$ and $k\geq 2$. Suppose that $H$ is a graph and $E(H, x) =E(G, x).$
Then $H$ is a $k$-regular graph of order $n$.
\end{teorem}

\nt We can determine some coefficients of $E(G, x)$ by having the degree
sequence of $G$. Also some information about the degree sequence of $G$ can be determined by the edge
cover polynomial of $G$.

\begin{teorem}\label{theorem3}{\rm \cite{saeid1}}
Let $G$ be a graph of order $n$ and size $m$ with no isolated vertex and
$E(G,x)=\sum_{ i=\rho(G)}^{m} e(G,i) x^{i}$, then the following hold:\\
\begin{enumerate}
\item[(i)] For every $i$,
\[
e(G, i) \geq {m \choose i} - \sum _{v \in V(G)}{ m-d(v) \choose i }.
\]
\item[(ii)] For every $i, i\geq m-2\delta +2$,
\[
e(G, i) = {m \choose i }- \sum _{v \in V(G)}{ m-d(v) \choose i}.
\]
\item[(iii)] For every $i, i\geq m-2\delta +2$,
\[
e(G, i+1) = \frac{m-i}{i+1}e(G,i)+\frac{1}{i+1} \sum _{v \in V(G)}{ m-d(v) \choose i}d(v).
\]
\item[(iv)] For every $k, 1\leq k\leq 2\delta -2$,
\[
a_k(G) = \frac{m-k+1}{k}e(G,m-k+1)- e(G,m-k)-\frac{1}{k}\sum _{j=1}^{k-1}{ m-j \choose m-k}j a_j(G).
\]
in particular if $\delta \geq 2$, then $a_{\delta}(G) ={m \choose {m-\delta}} - e(G,m - \delta)$.
\end{enumerate}
\end{teorem}

\begin{figure}[!h]
\hglue2.5cm
\includegraphics[width=11cm,height=5.1cm]{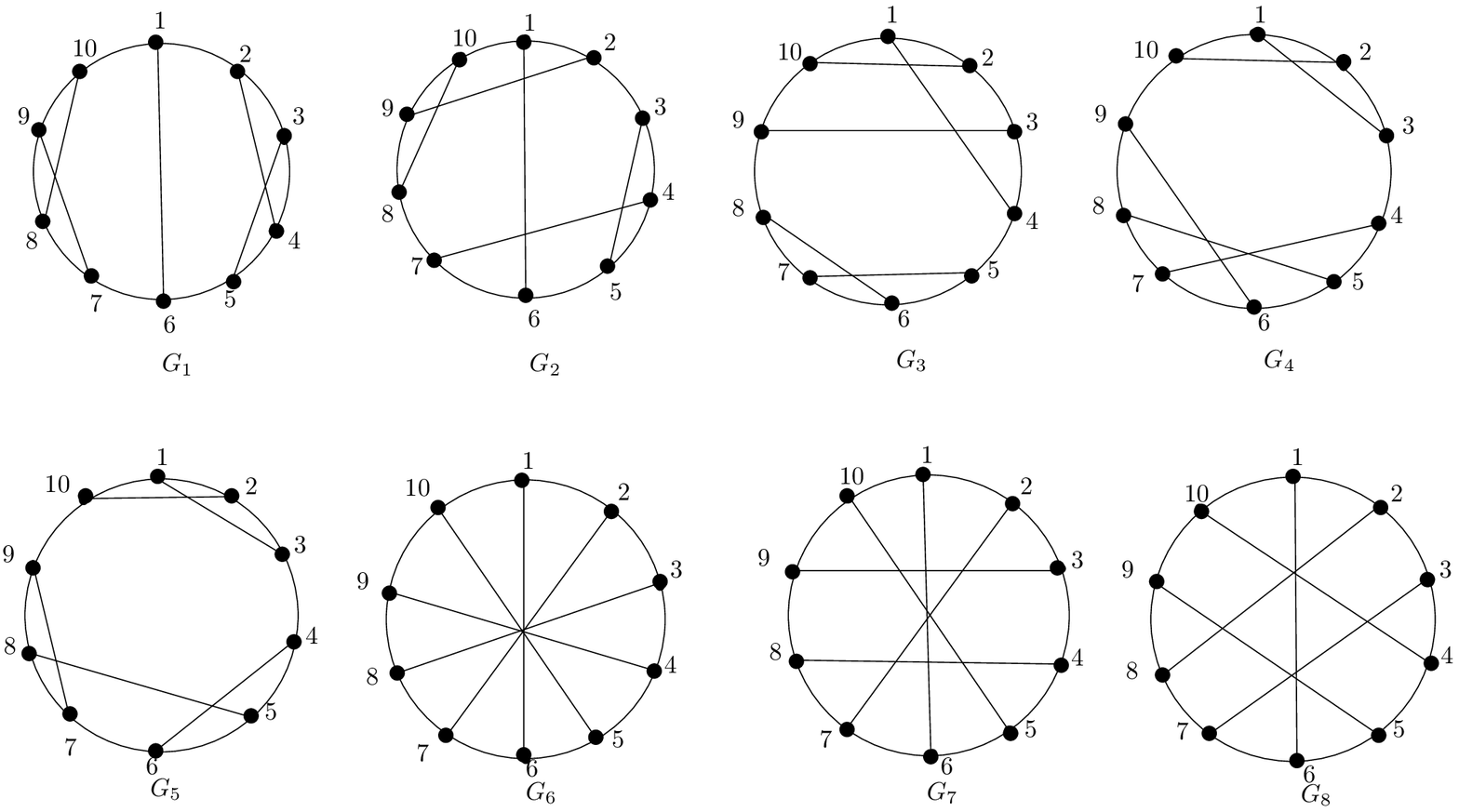}
\vglue5pt
\hglue2.5cm
\includegraphics[width=11cm,height=5.1cm]{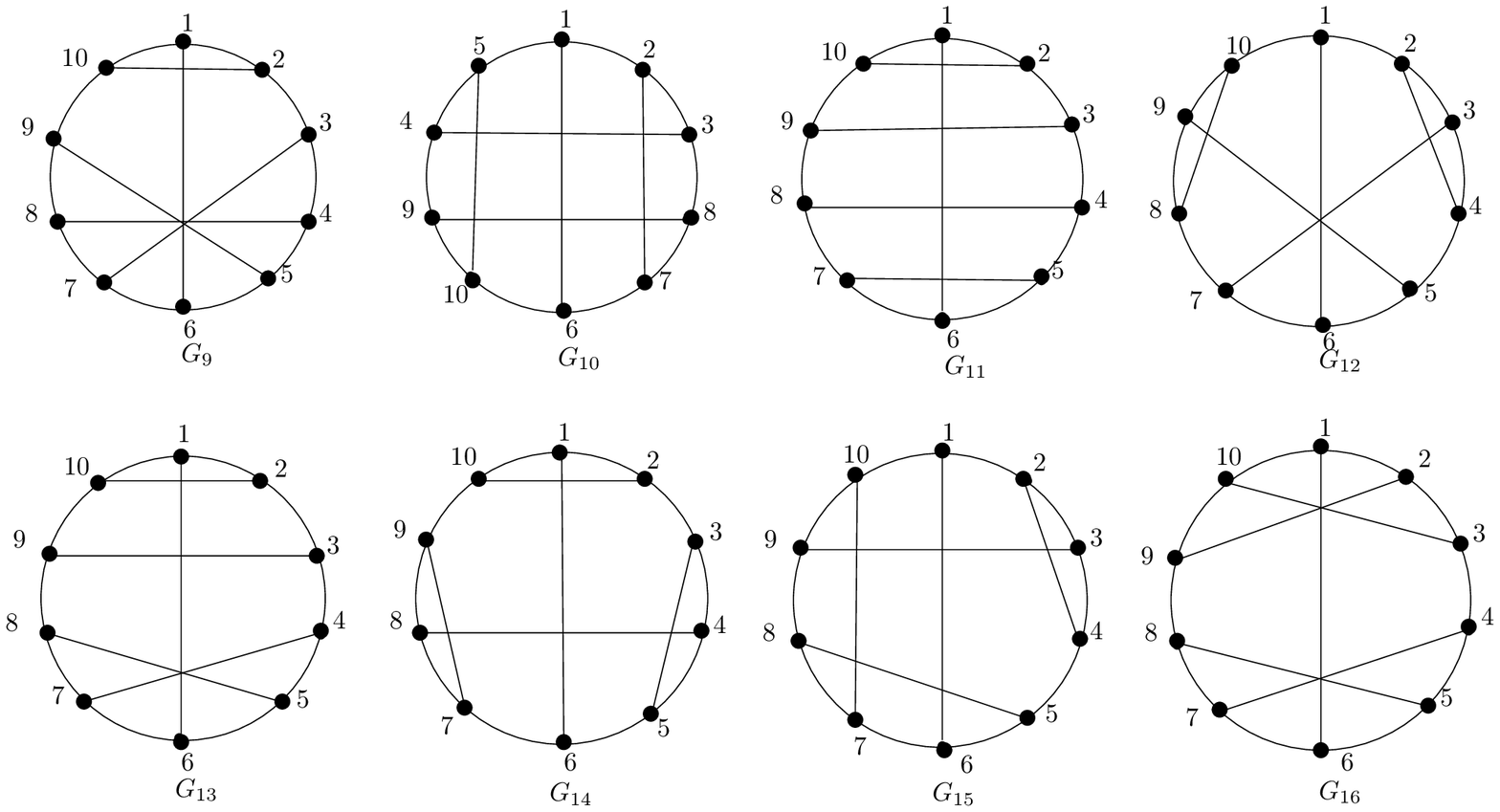}
\hglue2.5cm
\includegraphics[width=10.7cm,height=5cm]{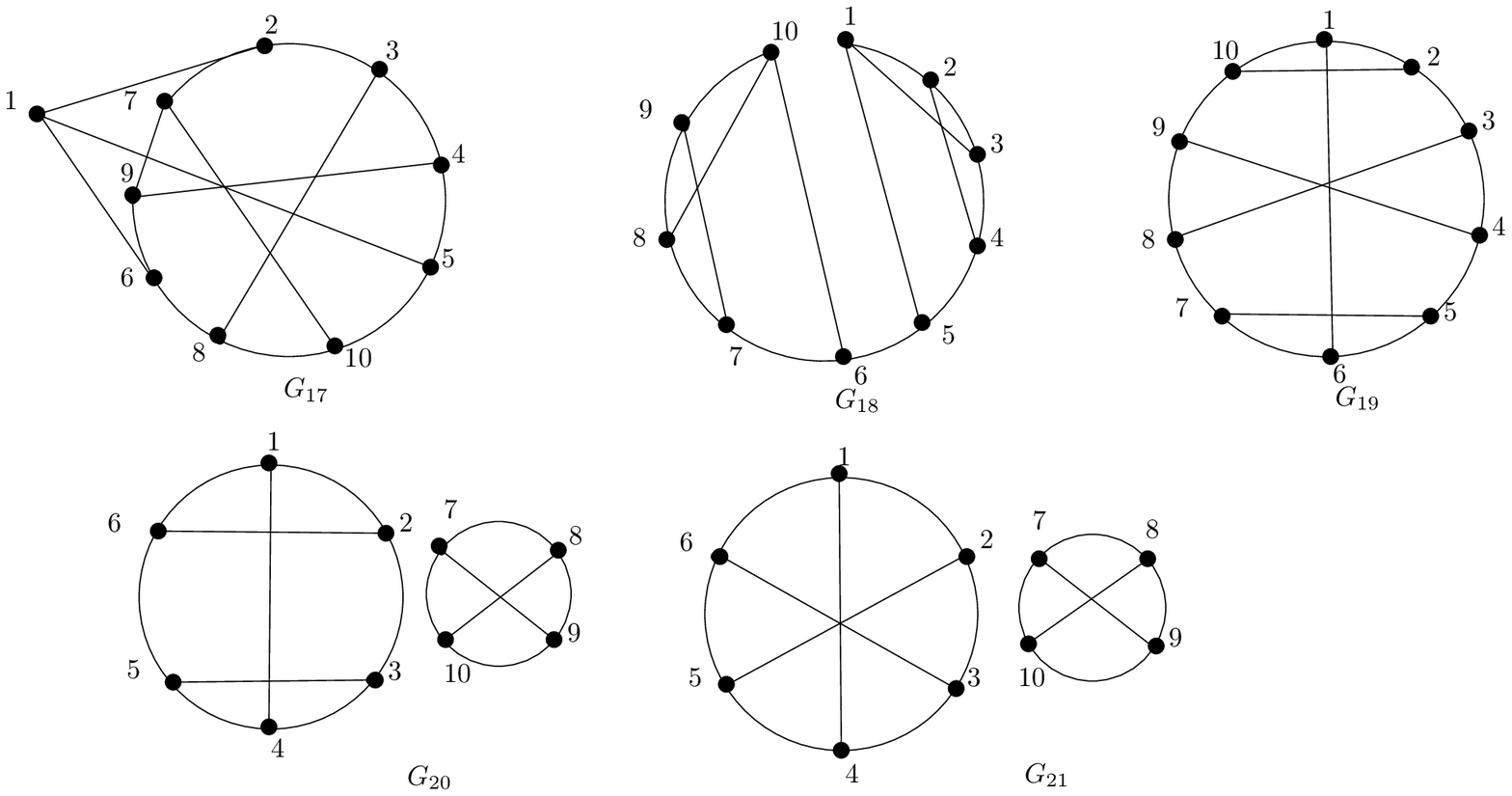}
\vglue-10pt \caption{\label{figure2} Cubic graphs of order 10.}
\end{figure}

\nt There are exactly $21$ cubic graphs of order $10$ given in
Figure~\ref{figure2} (see~\cite{reza}).
Clearly, edge covering number of a connected cubic
graph of order $10$ is $5$. Three are just two non-connected cubic graphs of order $10$.
Clearly, for these graphs, the edge covering number is also $5$. Note that the graph $G_{17}$ is
the Petersen graph. For the labeled graph $G_{17}$ in Figure~\ref{figure2}, we obtain all edge coverings
of size $5$ in the following lemma.

\begin{lema} \label{lemma2}
For the Petersen graph $P$, $e(P,5)=6$.
\end{lema}
\nt{\bf Proof.} We list all edge coverings of $P$ of
cardinality $5$, which are the $\rho$-sets of the labeled Petersen
graph (graph $G_{17}$) given in Figure~\ref{figure2}.
\\
$\mathscr{E}(P,5)=\Big\{\{\{1,2\},\{3,4\},\{5,10\},\{6,8\},\{7,9\}\},
\{\{1,2\},\{3,8\},\{4,5\},\{6,9\},\{7,10\}\},\\
\{\{1,5\},\{2,3\},\{4,9\},\{6,8\},\{7,10\}\},
\{\{1,5\},\{2,7\},\{3,4\},\{6,9\},\{8,10\}\},\\
\{\{1,6\},\{2,3\},\{4,5\},\{7,9\},\{8,10\}\},
\{\{1,6\},\{2,7\},\{3,8\},\{4,9\},\{5,10\}\}\Big\}.$
\nt Therefore, $e(P,5)=\big|\mathscr{E}(P,\rho)\big|=6$.\quad\qed

\nt We need the following theorem for finding the edge cover
polynomial of the Petersen graph.

\begin{lema}\label{lemma3}
Let $G$ be a cubic graph of order $10$. Then the following hold:
\begin{enumerate}
\item[(i)]
$e(G, i)={m \choose i}$, for $i=13,14,15$.
\item[(ii)]
$e(G, i)={m \choose i}-n{m-3 \choose i}$, for $i=11,12$.
\item[(iii)]
$e(G, 10)=
2358.$
\end{enumerate}
\end{lema}



\nt{\bf Proof.}
\begin{enumerate}
\item[(i)] It follows from  Lemma~\ref{lemma1} $(iii)$.

\item[(ii)]
It follows from   Theorem~\ref{theorem3} $(ii)$.

\item[(iii)] It is easy to see that $e(G, 10)={m \choose 10}-n{m-3 \choose 10}+15=
2358.$\quad\qed
\end{enumerate}

\begin{teorem} \label{theorem4}
The edge cover polynomial of the Petersen graph $P$ is:

$
E(P,x)=x^{15}+{15\choose 14}x^{14}+{15\choose 13}x^{13}+\Big[{15\choose 12} -10\Big]x^{12}+\Big[{15\choose 11}-120\Big]x^{11}+\Big[{15\choose 10}-645\Big]x^{10}
 +2985x^9+2400x^8+1101x^7+215x^6+6x^5.
$
\end{teorem}
\nt{\bf Proof.} The result follows from Lemmas~\ref{lemma2} and \ref{lemma3}.\quad\qed

\section{$\rho$ -sets of cubic graphs of order $10$.}

\nt In this section, we present all $\rho$-sets of connected cubic graphs
$G_1,G_2,\ldots,G_{18},G_{19}$ shown in
Figure~\ref{figure2}. The results here will be useful in studying the
$\mathscr{E}$-equivalence of these graphs in the last section.

\nt $\mathscr{E}(G_1,\rho)=\Big\{\{\{1,2\},\{3,4\},\{5,6\},\{7,8\},\{9,10\}\},
\{\{1,2\},\{3,4\},\{5,6\},\{7,9\},\{8,10\}\},\\
\{\{1,6\},\{2,3\},\{4,5\},\{7,8\},\{9,10\}\},
\{\{1,6\},\{2,3\},\{4,5\},\{7,9\},\{8,10\}\},\\
\{\{1,6\},\{2,4\},\{3,5\},\{7,8\},\{9,10\}\},
\{\{1,6\},\{2,4\},\{3,5\},\{7,9\},\{8,10\}\},\\
\{\{1,10\},\{2,3\},\{4,5\},\{6,7\},\{8,9\}\},
\{\{1,10\},\{2,4\},\{3,5\},\{6,7\},\{8,9\}\}\Big\}$.

\nt Therefore, $e(G_1,5)=\big|\mathscr{E}(G_1,\rho)\big|=8$.


\ms
\nt $\mathscr{E}(G_2,\rho)=\Big\{\{\{1,2\},\{3,4\},\{5,6\},\{7,8\},\{9,10\}\},
\{\{1,6\},\{2,3\},\{4,5\},\{7,8\},\{9,10\}\},\\
\{\{1,6\},\{2,9\},\{3,5\},\{4,7\},\{8,10\}\},
\{\{1,10\},\{2,3\},\{4,5\},\{6,7\},\{8,9\}\},\\
\{\{1,10\},\{2,3\},\{4,7\},\{5,6\},\{8,9\}\},
\{\{1,10\},\{2,9\},\{3,4\},\{5,6\},\{7,8\}\}\Big\}.$

\nt Therefore $e(G_2,5)=\big|\mathscr{E}(G_2,\rho)\big|=6$.

\ms

\nt $\mathscr{E}(G_3,\rho)=\Big\{\{\{1,2\},\{3,4\},\{5,6\},\{7,8\},\{9,10\}\},
\{\{1,2\},\{3,4\},\{5,7\},\{6,8\},\{9,10\}\},\\
\{\{1,4\},\{2,3\},\{5,6\},\{7,8\},\{9,10\}\},
\{\{1,4\},\{2,3\},\{5,7\},\{6,8\},\{9,10\}\},\\
\{\{1,4\},\{2,10\},\{3,9\},\{5,6\},\{7,8\}\},
\{\{1,4\},\{2,10\},\{3,9\},\{5,7\},\{6,8\}\},\\
\{\{1,10\},\{2,3\},\{4,5\},\{6,7\},\{8,9\}\}\Big\}.$
\nt Therefore $e(G_3,5)=\big|\mathscr{E}(G_3,\rho)\big|=7$.

\ms

\nt $\mathscr{E}(G_4,\rho)=\Big\{\{\{1,2\},\{3,4\},\{5,6\},\{7,8\},\{9,10\}\},
\{\{1,2\},\{3,4\},\{5,8\},\{6,7\},\{9,10\}\},\\
\{\{1,3\},\{2,10\},\{4,5\},\{6,7\},\{8,9\}\},
\{\{1,3\},\{2,10\},\{4,5\},\{6,9\},\{7,8\}\},\\
\{\{1,3\},\{2,10\},\{4,7\},\{5,6\},\{8,9\}\},
\{\{1,3\},\{2,10\},\{4,7\},\{5,8\},\{6,9\}\}, \\
\{\{1,10\},\{2,3\},\{4,5\},\{6,7\},\{8,9\}\},
\{\{1,10\},\{2,3\},\{4,5\},\{6,9\},\{7,8\}\},\\
\{\{1,10\},\{2,3\},\{4,7\},\{5,6\},\{8,9\}\},
\{\{1,10\},\{2,3\},\{4,7\},\{5,8\},\{6,9\}\}\Big\}.$

\nt Therefore $e(G_4,5)=\big|\mathscr{E}(G_4,\rho)\big|=10$.

\ms

\nt $\mathscr{E}(G_5,\rho)=\Big\{\{\{1,2\},\{3,4\},\{5,6\},\{7,8\},\{9,10\}\},
\{\{1,2\},\{3,4\},\{5,8\},\{6,7\},\{9,10\}\},\\
\{\{1,3\},\{2,10\},\{4,5\},\{6,7\},\{8,9\}\},
\{\{1,3\},\{2,10\},\{4,6\},\{5,8\},\{7,9\}\},\\
\{\{1,10\},\{2,3\},\{4,5\},\{6,7\},\{8,9\}\},
\{\{1,10\},\{2,3\},\{4,6\},\{5,8\},\{7,9\}\}\Big\}.$

\nt Therefore $e(G_5,5)=\big|\mathscr{E}(G_5,\rho)\big|=6$.

\ms

\nt $\mathscr{E}(G_6,\rho)=\Big\{\{\{1,2\},\{3,4\},\{5,6\},\{7,8\},\{9,10\}\},
\{\{1,2\},\{3,4\},\{5,10\},\{6,7\},\{8,9\}\},\\
\{\{1,2\},\{3,8\},\{4,5\},\{6,7\},\{9,10\}\},
\{\{1,2\},\{3,8\},\{4,9\},\{5,10\},\{6,7\}\},\\
\{\{1,6\},\{2,3\},\{4,5\},\{7,8\},\{9,10\}\},
\{\{1,6\},\{2,3\},\{4,9\},\{5,10\},\{7,8\}\},\\
\{\{1,6\},\{2,7\},\{3,4\},\{5,10\},\{8,9\}\},
\{\{1,6\},\{2,7\},\{3,8\},\{4,5\},\{9,10\}\},\\
\{\{1,6\},\{2,7\},\{3,8\},\{4,9\},\{5,10\}\},
\{\{1,10\},\{2,3\},\{4,5\},\{6,7\},\{8,9\}\}\\
\{\{1,10\},\{2,3\},\{4,9\},\{5,6\},\{7,8\}\},
\{\{1,10\},\{2,7\},\{3,4\},\{5,6\},\{8,9\}\},\\
\{\{1,10\},\{2,7\},\{3,8\},\{4,9\},\{5,6\}\}\Big\}.$
\nt Therefore $e(G_6,5)=\big|\mathscr{E}(G_6,\rho)\big|=13$.

\ms

\nt $\mathscr{E}(G_7,\rho)=\Big\{\{\{1,2\},\{3,4\},\{5,6\},\{7,8\},\{9,10\}\},
\{\{1,2\},\{3,4\},\{5,10\},\{6,7\},\{8,9\}\},\\
\{\{1,2\},\{3,9\},\{4,8\},\{5,10\},\{6,7\}\},
\{\{1,6\},\{2,3\},\{4,5\},\{7,8\},\{9,10\}\},\\
\{\{1,6\},\{2,7\},\{3,4\},\{5,10\},\{8,9\}\},
\{\{1,6\},\{2,7\},\{3,9\},\{4,8\},\{5,10\}\},\\
\{\{1,10\},\{2,3\},\{4,5\},\{6,7\},\{8,9\}\},
\{\{1,10\},\{2,7\},\{3,4\},\{5,6\},\{8,9\}\},\\
\{\{1,10\},\{2,7\},\{3,9\},\{4,8\},\{5,6\}\}\Big\}.$
\nt Therefore $e(G_7,5)=\big|\mathscr{E}(G_7,\rho)\big|=9$

\ms

\nt $\mathscr{E}(G_8,\rho)=\Big\{\{\{1,2\},\{3,4\},\{5,6\},\{7,8\},\{9,10\}\},
\{\{1,2\},\{3,7\},\{4,10\},\{5,6\},\{8,9\}\},\\
\{\{1,6\},\{2,3\},\{4,5\},\{7,8\},\{9,10\}\},
\{\{1,6\},\{2,3\},\{4,10\},\{5,9\},\{7,8\}\},\\
\{\{1,6\},\{2,8\},\{3,7\},\{4,5\},\{9,10\}\},
\{\{1,6\},\{2,8\},\{3,7\},\{4,10\},\{5,9\}\},\\
\{\{1,10\},\{2,3\},\{4,5\},\{6,7\},\{8,9\}\},
\{\{1,10\},\{2,8\},\{3,4\},\{5,9\},\{6,7\}\}\Big\}$.

\nt Therefore $e(G_8,5)=\big|\mathscr{E}(G_8,\rho)\big|=8.$

\ms

\nt $\mathscr{E}(G_9,\rho)=\Big\{\{\{1,2\},\{3,4\},\{5,6\},\{7,8\},\{9,10\}\},
\{\{1,2\},\{3,7\},\{4,8\},\{5,6\},\{9,10\}\},\\
\{\{1,6\},\{2,3\},\{4,5\},\{7,8\},\{9,10\}\},
\{\{1,6\},\{2,10\},\{3,4\},\{5,9\},\{7,8\}\},\\
\{\{1,6\},\{2,10\},\{3,7\},\{4,5\},\{8,9\}\},
\{\{1,6\},\{2,10\},\{3,7\},\{4,8\},\{5,9\}\},\\
\{\{1,10\},\{2,3\},\{4,5\},\{6,7\},\{8,9\}\},
\{\{1,10\},\{2,3\},\{4,8\},\{5,9\},\{6,7\}\}\Big\}$.

\nt Therefore $e(G_9,5)=\big|\mathscr{E}(G_9,\rho)\big|=8.$

\ms

\nt $\mathscr{E}(G_{10},\rho)=\Big\{\{\{1,2\},\{3,4\},\{5,10\},\{6,7\},\{8,9\}\},
\{\{1,2\},\{3,8\},\{4,5\},\{6,7\},\{9,10\}\},\\
\{\{1,2\},\{3,8\},\{4,9\},\{5,10\},\{6,7\}\},
\{\{1,5\},\{2,3\},\{4,9\},\{6,10\},\{7,8\}\},\\
\{\{1,5\},\{2,7\},\{3,4\},\{6,10\},\{8,9\}\},
\{\{1,5\},\{2,7\},\{3,8\},\{4,9\},\{6,10\}\},\\
\{\{1,6\},\{2,3\},\{4,5\},\{7,8\},\{9,10\}\},
\{\{1,6\},\{2,3\},\{4,9\},\{5,10\},\{7,8\}\},\\
\{\{1,6\},\{2,7\},\{3,4\},\{5,10\},\{8,9\}\},
\{\{1,6\},\{2,7\},\{3,8\},\{4,5\},\{9,10\}\},\\
\{\{1,6\},\{2,7\},\{3,8\},\{4,9\},\{5,10\}\}\Big\}$.
\nt Therefore $e(G_{10},5)=\big|\mathscr{E}(G_{10},\rho)\big|=11.$

\ms

\nt $\mathscr{E}(G_{11},\rho)=\Big\{\{\{1,2\},\{3,4\},\{5,6\},\{7,8\},\{9,10\}\},
\{\{1,6\},\{2,3\},\{4,5\},\{7,8\},\{9,10\}\},\\
\{\{1,6\},\{2,3\},\{4,8\},\{5,7\},\{9,10\}\},
\{\{1,6\},\{2,10\},\{3,4\},\{5,7\},\{8,9\}\},\\
\{\{1,6\},\{2,10\},\{3,9\},\{4,5\},\{7,8\}\},
\{\{1,6\},\{2,10\},\{3,9\},\{4,8\},\{5,7\}\},\\
\{\{1,10\},\{2,3\},\{4,5\},\{6,7\},\{8,9\}\}\Big\}$.
\nt Therefore $e(G_{11},5)=\big|\mathscr{E}(G_{11},\rho)\big|=7$.

\ms

\nt $\mathscr{E}(G_{12},\rho)=\Big\{\{\{1,2\},\{3,4\},\{5,6\},\{7,8\},\{9,10\}\},
\{\{1,2\},\{3,4\},\{5,9\},\{6,7\},\{8,10\}\},\\
\{\{1,6\},\{2,3\},\{4,5\},\{7,8\},\{9,10\}\},
\{\{1,6\},\{2,4\},\{3,7\},\{5,9\},\{8,10\}\},\\
\{\{1,10\},\{2,3\},\{4,5\},\{6,7\},\{8,9\}\},
\{\{1,10\},\{2,4\},\{3,7\},\{5,6\},\{8,9\}\}\Big\}.$

Therefore $e(G_{12},5)=\big|\mathscr{E}(G_{12},\rho)\big|=6$.

\ms

\nt $\mathscr{E}(G_{13},\rho)=\Big\{\{\{1,2\},\{3,4\},\{5,6\},\{7,8\},\{9,10\}\},
\{\{1,2\},\{3,4\},\{5,8\},\{6,7\},\{9,10\}\},\\
\{\{1,6\},\{2,3\},\{4,5\},\{7,8\},\{9,10\}\},
\{\{1,6\},\{2,3\},\{4,7\},\{5,8\},\{9,10\}\},\\
\{\{1,6\},\{2,10\},\{3,9\},\{4,5\},\{7,8\}\},
\{\{1,6\},\{2,10\},\{3,9\},\{4,7\},\{5,8\}\},\\
\{\{1,10\},\{2,3\},\{4,5\},\{6,7\},\{8,9\}\},
\{\{1,10\},\{2,3\},\{4,7\},\{5,6\},\{8,9\}\}\Big\}$.

\nt Therefore $e(G_{13},5)=\big|\mathscr{E}(G_{13},\rho)\big|=8$.

\ms

\nt $\mathscr{E}(G_{14},\rho)=\Big\{\{\{1,2\},\{3,4\},\{5,6\},\{7,8\},\{9,10\}\},
\{\{1,2\},\{3,5\},\{4,8\},\{6,7\},\{9,10\}\},\\
\{\{1,6\},\{2,3\},\{4,5\},\{7,8\},\{9,10\}\},
\{\{1,6\},\{2,10\},\{3,5\},\{4,8\},\{7,9\}\},\\
\{\{1,10\},\{2,3\},\{4,5\},\{6,7\},\{8,9\}\},
\{\{1,10\},\{2,3\},\{4,8\},\{5,6\},\{7,9\}\}\Big\}$.

\nt Therefore $e(G_{14},5)=\big|\mathscr{E}(G_{14},\rho)\big|=6$.

\ms

\nt $\mathscr{E}(G_{15},\rho)=\Big\{\{\{1,2\},\{3,4\},\{5,6\},\{7,8\},\{9,10\}\},
\{\{1,2\},\{3,4\},\{5,6\},\{7,10\},\{8,9\}\},\\
\{\{1,2\},\{3,4\},\{5,8\},\{6,7\},\{9,10\}\},
\{\{1,6\},\{2,3\},\{4,5\},\{7,8\},\{9,10\}\},\\
\{\{1,6\},\{2,3\},\{4,5\},\{7,10\},\{8,9\}\},
\{\{1,6\},\{2,4\},\{3,9\},\{5,8\},\{7,10\}\},\\
\{\{1,10\},\{2,3\},\{4,5\},\{6,7\},\{8,9\}\},
\{\{1,10\},\{2,4\},\{3,9\},\{5,6\},\{7,8\}\},\\
\{\{1,10\},\{2,4\},\{3,9\},\{5,8\},\{6,7\}\}\Big\}$.
\nt Therefore $e(G_{15},5)=\big|\mathscr{E}(G_{15},\rho)\big|=9$.

\ms

\nt $\mathscr{E}(G_{16},\rho)=\Big\{ \{\{1,2\},\{3,4\},\{5,6\},\{7,8\},\{9,10\}\},
\{\{1,2\},\{3,4\},\{5,8\},\{6,7\},\{9,10\}\},\\
\{\{1,2\},\{3,10\},\{4,5\},\{6,7\},\{8,9\}\},
\{\{1,2\},\{3,10\},\{4,7\},\{5,6\},\{8,9\}\},\\
\{\{1,6\},\{2,3\},\{4,5\},\{7,8\},\{9,10\}\},
\{\{1,6\},\{2,3\},\{4,7\},\{5,8\},\{9,10\}\},\\
\{\{1,6\},\{2,9\},\{3,10\},\{4,5\},\{7,8\}\},
\{\{1,6\},\{2,9\},\{3,10\},\{4,7\},\{5,8\}\},\\
\{\{1,10\},\{2,3\},\{4,5\},\{6,7\},\{8,9\}\},
\{\{1,10\},\{2,3\},\{4,7\},\{5,6\},\{8,9\}\},\\
\{\{1,10\},\{2,9\},\{3,4\},\{5,6\},\{7,8\}\},
\{\{1,10\},\{2,9\},\{3,4\},\{5,8\},\{6,7\}\}\Big\}$.

\nt Therefore $e(G_{16},5)=\big|\mathscr{E}(G_{16},\rho)\big|=12$.

\ms

\nt $\mathscr{E}(G_{17},\rho)=\Big\{\{\{1,2\},\{3,4\},\{5,10\},\{6,8\},\{7,9\}\},
\{\{1,2\},\{3,8\},\{4,5\},\{6,9\},\{7,10\}\},\\
\{\{1,5\},\{2,3\},\{4,9\},\{6,8\},\{7,10\}\},
\{\{1,5\},\{2,7\},\{3,4\},\{6,9\},\{8,10\}\},\\
\{\{1,6\},\{2,3\},\{4,5\},\{7,9\},\{8,10\}\},
\{\{1,6\},\{2,7\},\{3,8\},\{4,9\},\{5,10\}\}\Big\}.$

\nt Therefore, $e(G_{17},5)=\big|\mathscr{E}(G_{17},\rho)\big|=6$.

\ms

\nt $\mathscr{E}(G_{18},\rho)=\Big\{\{\{1,2\},\{3,4\},\{5,6\},\{7,8\},\{9,10\}\},
\{\{1,2\},\{3,4\},\{5,6\},\{7,9\},\{8,10\}\},\\
\{\{1,3\},\{2,4\},\{5,6\},\{7,8\},\{9,10\}\},
\{\{1,3\},\{2,4\},\{5,6\},\{7,9\},\{8,10\}\}\Big\}.$

\nt Therefore, $e(G_{18},5)=\big|\mathscr{E}(G_{18},\rho)\big|=4$.

\ms

\nt $\mathscr{E}(G_{19},\rho)=\Big\{\{\{1,2\},\{3,4\},\{5,6\},\{7,8\},\{9,10\}\},
\{\{1,2\},\{3,8\},\{4,5\},\{6,7\},\{9,10\}\},\\
\{\{1,6\},\{2,3\},\{4,5\},\{7,8\},\{9,10\}\},
\{\{1,6\},\{2,10\},\{3,4\},\{5,7\},\{8,9\}\},\\
\{\{1,6\},\{2,10\},\{3,8\},\{4,9\},\{5,7\}\},
\{\{1,10\},\{2,3\},\{4,5\},\{6,7\},\{8,9\}\},\\
\{\{1,10\},\{2,3\},\{4,9\},\{5,6\},\{7,8\}\}\Big\}$.
\nt Therefore $e(G_{19},5)=\big|\mathscr{E}(G_{19},\rho)\big|=7$.

\section{${\mathscr{E}}$-Equivalence classes of all cubic graphs of order $10$}

\nt In this section we study the ${\mathscr{E}}$-equivalence classes
of all cubic graphs of order $10$ and we show
that all cubic graphs of order $10$ are determined uniquely by their edge cover polynomials.

\begin{teorem}
Let $G$ be a cubic graphs of order $10$ and $H$ be a graph of order $n$. If  $E(H,x) = E(G, x),$ then $H=G$.
\end{teorem}

\nt{\bf Proof.} By Theorem \ref{theorem2}, we find that $H$ is a $3$-regular graph of order $10$. Since there are
 exactly $21$ cubic graphs of order $10$ (Figure \ref{figure2}), by comparing the edge cover polynomial of these graphs,
we conclude  that none of these polynomials are equal. Therefore $H=G$.\quad\qed

\bigskip
\nt The following procedure allows to compute the edge cover polynomial of a graph $G$ with order $n$ and size $m$. Let $W = \{S_1, \cdots , S_{2^m-1}\}$ be the collection of nonempty subsets of $E(G)$ .
The algorithm starts with $E(G,x) = 0$ and continues with the following steps, for $1 \leq j \leq 2^m - 1.$
\begin{itemize}
\item[1.] If $\bigcup_{i=1}^{|S_j|} {S_j}_i =\{1, \cdots ,n\}$, then go to step $2$, else replace $j$ by $j + 1$ and apply the algorithm again.
\item[2.] Add one term $x^{|S_j|}$ to $E(G, x)$.
\item[3.] Replace $j$ by $j + 1$ and apply the algorithm again.
\end{itemize}

\nt We show the number of edge covering  of the cubic graphs of order $10$  in Table 1.

\begin{center}
\begin{footnotesize}
\small
\begin{tabular}{r||cccccccccccc}
$j$&$5$&$6$&$7$&$8$&$9$&$10$&$11$&$12$&$13$&$14$&$15$&\\ [0.1ex]
\hline
$e(G_1,j)$&8&195&1055&2378&2981&2358&1245&445&105&15&1&\\
$e(G_2,j)$&6&205&1073&2388&2983&2358&1245&445&105&15&1&\\
$e(G_3,j)$&7&202&1065&2383&2982&2358&1245&445&105&15&1&\\
$e(G_4,j)$&10&218&1080&2389&2983&2358&1245&445&105&15&1&\\
$e(G_5,j)$&6&192&1052&2377&2981&2358&1245&445&105&15&1&\\
$e(G_6,j)$&13&230&1100&2400&2985&2358&1245&445&105&15&1&\\
$e(G_7,j)$&9&224&1098&2400&2985&2358&1245&445&105&15&1&\\
$e(G_8,j)$&8&221&1097&2400&2985&2358&1245&445&105&15&1&\\
$e(G_9,j)$&8&214&1085&2394&2984&2358&1245&445&105&15&1&\\
$e(G_{10},j)$&11&230&1100&2400&2985&2358&1245&445&105&15&1&\\
$e(G_{11},j)$&7&208&1074&2388&2983&2358&1245&445&105&15&1&\\
$e(G_{12},j)$&6&201&1071&2388&2983&2358&1245&445&105&15&1&\\
$e(G_{13},j)$&8&218&1087&2394&2984&2358&1245&445&105&15&1&\\
$e(G_{14},j)$&6&194&1059&2382&2982&2358&1245&445&105&15&1&\\
$e(G_{15},j)$&9&217&1086&2394&2984&2358&1245&445&105&15&1&\\
$e(G_{16},j)$&12&231&1101&2400&2985&2358&1245&445&105&15&1&\\
$e(G_{17},j)$&6&215&1095&2400&2985&2358&1245&445&105&15&1&\\
$e(G_{18},j)$&4&197&1057&2378&2981&2358&1245&445&105&15&1&\\
$e(G_{19},j)$&7&204&1072&2388&2983&2358&1245&445&105&15&1&\\
$e(G_{20},j)$&12&199&1050&2373&2980&2358&1245&445&105&15&1&\\
$e(G_{21},j)$&18&231&1080&2385&2982&2358&1245&445&105&15&1&
\end{tabular}
\end{footnotesize}
\end{center}
\begin{center}
\nt{Table 1.} $e(G_{n},j)$, the number of edge coverings of $G_n$ with cardinality $j$.
\end{center}
\nt Regarding to this table we see that
all cubic graphs of order $10$ are determined uniquely by their edge cover polynomial. So we have the following corollary.

\begin{korolari}
All cubic graphs of order $10$ are determined uniquely by their edge cover polynomials.
\end{korolari}

\end{document}